\documentclass[12pt]{amsart}
\usepackage{a4wide}
\usepackage{amssymb}
\usepackage{amsthm}
\usepackage{amsmath}
\usepackage{amscd}
\usepackage{verbatim}

\numberwithin{equation}{section}

\theoremstyle{plain}
\newtheorem{theorem}{Theorem}[section]

\newtheorem{proposition}[theorem]{Proposition}

\theoremstyle{definition}

\newtheorem{remark}[theorem]{Remark}

\theoremstyle{remark}

\newcommand{\R}{\mathbb{R}}
\newcommand{\Q}{\mathbb{Q}}
\newcommand{\Z}{\mathbb{Z}}
\newcommand{\N}{\mathbb{N}}
\newcommand{\C}{\mathbb{C}}

\renewcommand{\H}{\mathbb{H}}

\newcommand{\D}{\mathbb{D}}

\newcommand{\leg}[2]{\left( \frac{#1}{#2} \right)}

\newcommand{\eps}{\varepsilon}

\begin{document}

\title[Sign changes of coefficients modular forms]{Sign changes of coefficients of half integral weight modular forms}

\author{Jan Hendrik Bruinier and Winfried Kohnen}

\address{Technische Universit\"at Darmstadt,
Fachbereich Mathematik,
Schlossgartenstrasse 7,
D--64289 Darmstadt,
Germany}
\email{bruinier@mathematik.tu-darmstadt.de }

\address{Mathematisches Institut, Universit\"at Heidelberg, Im Neuenheimer Feld 288, D-69120 Heidelberg, Germany}
\email{winfried@mathi.uni-heidelberg.de}


\date{September 3, 2007}

\begin{abstract}
For a half integral weight modular form $f$ we study the signs of the 
Fourier coefficients $a(n)$.
If $f$ is a Hecke eigenform of level $ N$ with real Nebentypus character, 
and $t$ is a fixed square-free  positive integer with 
$a(t)\neq 0$, we show that for all but finitely many primes $p$ the sequence 
$(a(tp^{2m}))_{m}$ has infinitely many signs changes. Moreover, we prove
similar (partly conditional) results for arbitrary cusp forms $f$ 
which are not necessarily Hecke eigenforms.  
\end{abstract} 

\maketitle

\section{Introduction}

Let $f$ be a non-zero elliptic cusp form of positive real weight
$\kappa$, with multiplier $v$ and with real Fourier coefficients
$a(n)$ for  $n\in {\N}$. Then under quite general conditions, using
the theory of $L$-functions, it was shown in \cite{KKP} that the
sequence $(a(n))_{n\in {\N}}$ has infinitely many sign changes,
i.e., there are infinitely many $n$ such that $a(n)>0$ and there are
infinitely many $n$ such that $a(n)<0$.

This is particularly interesting when $\kappa$ is an integer and $f$
is a Hecke eigenform of level $N$, and so the $a(n)$ are
proportional to the Hecke eigenvalues. For recent work in this
direction we refer to e.g. \cite{KoSe}, \cite{IKS}, \cite{KSL}.

In the present note we shall consider the case of half-integral
weight $\kappa=k+1/2$, $k\in {\N}$, and level $N$ divisible by $4$.
Note that this case is distinguished through the celebrated works of
Shimura \cite{Sh} and Waldspurger \cite{Wa} in the following way.
First, for each square-free positive integer $t$, there exists a
linear lifting from weight $k+1/2$ to even integral weight $2k$
determined by the coefficients $a(tn^2)$ (where $n\in {\N}$), see
\cite{Ni}, \cite{Sh}. In particular, through these liftings, the
theory of Hecke eigenvalues is the same as that in the integral
weight case. Secondly, if $f$ is a Hecke eigenform, then the {\em
squares} $a(t)^2$ are essentially proportional to the central
critical values of the Hecke $L$-function of $F$ twisted with the
quadratic character $\chi_{t,N}= (\frac{(-1)^kN^2t}{\cdot})$, see
\cite{Wa}. Here $F$ is a Hecke eigenform of weight $2k$
corresponding to $f$ under the Shimura correspondence.

These facts motivate the following questions. First, it is  natural
to ask for sign changes of the sequence $(a(tn^2))_{n\in {\N}}$
where $t$ is a fixed positive square-free integer. 
We start with a conditional result here, 
namely if the Dirichlet $L$-function
associated to $\chi_{t,N}$ has no zeros in the interval $(0,1)$
(Chowla's conjecture), then the sequence $(a(tn^2))_{n\in {\N}}$
--if not identically zero-- changes sign infinitely often (Theorem
\ref{thm2}). Note that by work of Conrey and Soundararajan
\cite{CS}, Chowla's conjecture is true for a positive proportion of
positive square-free integers $t$. 
If $f$ is a Hecke eigenform, we can in fact prove an unconditional 
and much better result on the sign changes of the sequence 
$(a(tp^{2m}))_{m\in {\N}}$ where $p$ is a prime not dividing $N$
(see Theorem \ref{thm2.5}).

Secondly, one may ask for sign changes of the sequence $( a(t) )_t$,
where $t$ runs through the square-free integers only. This question
is more difficult to treat. Numerical calculations
seem to suggest not only that there are infinitely
many sign changes, but also that ``half'' of these coefficients are
positive and ``half'' of them are negative.

It seems quite difficult to prove any general theorem here, and we
can only prove a result that seems to point into the right direction (see
Theorem \ref{thm1}): Under the (clearly necessary) assumption for
$k=1$ that $f$ is contained in the orthogonal complement of the
space of unary theta functions, there exist infinitely many positive
square-free integers $t$ and for each such $t$ a natural number
$n_t$, such that the sequence $( a(tn_t^2) )_t$ has infinitely many
sign changes. Note that we do not require $f$ to be an eigenform. We
in fact prove a slightly stronger result (Theorem \ref{cor1}).

As an immediate application, in the integral weight as well as in
the half-integral weight case, we may consider representation
numbers of quadratic forms. Let $Q$ be a positive definite integral
quadratic form, and for a positive integer $n$ let $r_Q(n)$ be the
number of integral representations of $n$ by $Q$. Then the
associated theta series is the sum of a modular form lying in the
space of Eisenstein series and a corresponding cusp form. An
infinity of sign changes of the coefficients of the latter (the
``error term'' for $r_Q(n)$) means that $r_Q(n)$ for infinitely many
$n$ is larger (respectively less) than the corresponding Eisenstein
coefficient (the ``main term'' for $r_Q(n)$).

Exact statements of our results are given in Section 2, while
Section 3 contains their proofs. These are based on the existence of
the Shimura lifts, the theory of $L$-functions, and on results on
quadratic twists proved in \cite{Br}. In Section 4 some numerical
examples are given.


\section{Notation and statement of results}

We denote by $\N$ the set of positive integers. The set of
square-free positive integers is denoted by $\D$. Throughout we
write $q=e^{2\pi i
  z}$ for $z$ in the upper complex half plane $\H$.

Let $k$ be a positive integer. Let $N$ be a positive integer
divisible by $4$, and let $\chi$ be a Dirichlet character modulo
$N$.  We write $\chi^*$ for the Dirichlet character modulo $N$ given
by $\chi^*(a)=\leg{-4}{a}^k\chi(a)$.  Moreover, we write
$S_{k+1/2}(N,\chi)$ for the space of cusp forms of weight $k+1/2$
for the group $\Gamma_0(N)$ with character $\chi$ in the sense of
Shimura \cite{Sh}.

If $m$ and $r$ are positive integers and $\psi$ is an odd primitive
Dirichlet character modulo $r$, then the unary theta function
\[
\theta_{\psi,m}(z)=\sum_{n\in \Z} \psi(n) n q^{m n^2}
\]
belongs to $S_{3/2}\left(N, \leg{-4m}{\cdot}\psi\right)$ for all $N$
divisible by $4r^2m$, cf.~\cite{Sh}. Let $S_{3/2}^*(N,\chi)$ be the
orthogonal complement with respect to the Petersson scalar product
of the subspace of $S_{3/2}(N,\chi)$ spanned by such theta series.
For $k\geq 2$ we simply put $S_{k+1/2}^*(N,\chi)=S_{k+1/2}(N,\chi)$.
It is well known that the Shimura lift maps $S^*_{k+1/2}(N,\chi)$ to
the space $S_{2k}(N/2,\chi^2)$ of cusp forms of integral weight $2k$
for $\Gamma_0(N/2)$ with character $\chi^2$.

Throughout this section, let $f=\sum_{n=1}^\infty a(n) q^n\in
S^*_{k+1/2}(N,\chi)$ be a non-zero cup form with Fourier
coefficients $a(n)\in \R$. In our first result we consider the
coefficients $a(t n^2)$ for fixed $t\in \D$ and varying $n$.

\begin{theorem}
\label{thm2}
Let $t\in \D$ such that $a(t)\neq 0$, and write $\chi_{t,N}$ for the
quadratic character $\chi_{t,N}= (\frac{(-1)^kN^2t}{\cdot})$. Assume
that the Dirichlet $L$-function $L(s,\chi_{t,N})$ has no zeros in
the interval $(0,1)$. Then the sequence $(a(tn^2))_{n\in \N}$ has
infinitely many sign changes.
%
\end{theorem}

For Hecke eigenforms, we prove the following unconditional result.

\begin{theorem}
\label{thm2.5}
Suppose that the character $\chi$ of $f$ is real, and suppose 
that $f$ is an eigenform of all 
Hecke operators $T(p^2)$ with corresponding eigenvalues 
$\lambda_p$ for $p$ coprime to $N$.
Let $t\in \D$ such that 
$a(t)\neq 0$. Then for all but finitely many  primes $p$ coprime to $N$
the sequence $(a(tp^{2m}))_{m\in {\N}}$ has infinitely 
many sign changes. 
\end{theorem}

\begin{remark}
\label{rem2.5}
Let $K_f$ be the number field generated by the Hecke eigenvalues $\lambda_p$ of $f$.
The number of exceptional primes in Theorem \ref{thm2.5} is bounded by $r$ 
where $2^r$ is the highest power of $2$ dividing the degree of $K_f$ over $\Q$.
\end{remark}

Next, we consider the coefficients $a(t n^2)$ for varying $t\in \D$.

\begin{theorem}
\label{thm1}
For every $t\in \D$ there is an $n_t\in \N$ such that the sequence
$(a(t n_t^2))_{t\in \D}$ has infinitely many sign changes.
\end{theorem}

In the special case when $f$ is a Hecke eigenform, Theorem \ref{thm1} 
is an easy consequence of 
Theorem \ref{thm2.5}. However, we do not assume this in  Theorem \ref{thm1}.
Notice that the statement is obviously wrong for the theta functions
$\theta_{\psi,m}$. We shall also prove the following slightly
stronger statement.

\begin{theorem}
\label{cor1}
Let $p_1,\dots,p_r$ be distinct primes not dividing $N$ and let
$\eps_1,\dots,\eps_r\in \{\pm 1\}$. Write $\D'$ for the set of $t\in
\D$ satisfying $(\frac{t}{p_j})=\eps_j$ for $j=1,\dots,r$.
Then for every $t\in \D'$ there is an $n_t\in \N$ such that the
sequence $(a(t n_t^2))_{t\in \D'}$ has infinitely many sign changes.
%
\end{theorem}

\section{Proofs}

Here we prove the results of the previous section.

\begin{proof}[Proof of Theorem \ref{thm2}.]
Put
\begin{align}
\label{w1} A(n):= \sum_{d|n}\chi_{t,N}(d)d^{k-1}a(\frac{n^2}{d^2}t).
\end{align}
According to \cite{Sh,Ni}, the series
\[
F(z):=\sum_{n\geq 1}A(n)q^n
\]
is in $S_{2k}(N/2,\chi^2)$ and is non-zero due to our assumption
$a(t)\neq 0$. Note that \eqref{w1} is equivalent to the Dirichlet
series identity
\begin{align}
\label{w2} \sum_{n\geq
1}a(tn^2)n^{-s}=\frac{1}{L(s-k+1,\chi_{t,N})}\cdot L(F,s)
\end{align}
in the range of absolute convergence, where $L(F,s)$ is the Hecke
$L$-function attached to $F$.

Now suppose that $a(tn^2)\geq 0$ for all but finitely many $n$. Then
by a classical theorem of Landau, either the Dirichlet series on the
left hand side of \eqref{w2} has a singularity at the real point of
its line of convergence or must converge everywhere.

By our hypothesis, $L(s,\chi_{t,N})$ has no real zeros for
$\Re(s)>0$. Hence the series on the left hand side of
\eqref{w2} converges for $\Re(s)>k-1$. In particular, we have
\[
a(tn^2)\ll_{\epsilon} n^{k-1+\epsilon}\qquad (\epsilon>0).
\]
From \eqref{w1} we therefore deduce that
\[
A(n)\ll_\epsilon
\sum_{d|n}d^{k-1}\left(\frac{n}{d}\right)^{k-1+\epsilon}\ll_\epsilon
n^{k-1+2\epsilon} \qquad (\epsilon>0).
\]
Consequently, the Rankin-Selberg Dirichlet series
\[
R_F(s)=\sum_{n\geq 1}A(n)^2n^{-s}
\]
must be convergent for $\Re(s)>2k-1$. However, it is well-known that
the latter has a pole at $s=2k$ with residue $c_k \|F\|^2$, where
$c_k>0$ is a constant depending only on $k$, and $\|F\|^2$ is the
square of the Petersson norm of $F$. Since $F\neq 0$, we obtain a
contradiction. This proves the claim.
\end{proof}

\begin{proof}[Proof of Theorem \ref{thm2.5} and Remark \ref{rem2.5}.]
We use the same notation as in the proof of Theorem \ref{thm2}.
Since $f$ is an eigenfunction of $T(p^2)$, the function $F$ is an 
eigenfunction under the usual Hecke operator $T(p)$ with eigenvalue 
$\lambda_p$. Since $\chi^2=1$,  the eigenvalue $\lambda_p$ is real.
One has
\begin{align}
\label{w3.3}
\sum_{m\geq 1}a(tp^{2m})p^{-ms} = 
a(t)\frac{1-\chi_{t,N}(p)p^{k-1-s}}{1-\lambda_p p^{-s}+p^{2k-1-2s}}
\end{align}
for $\Re(s)$ sufficiently large, which is the local variant of \eqref{w2}.

The denominator of the right-hand side of \eqref{w3.3} factors as
\[
1-\lambda_pp^{-s}+p^{2k-1-2s}=(1-\alpha_pp^{-s})(1-\beta_pp^{-s})
\]
where $\alpha_p+\beta_p=\lambda_p$ and $\alpha_p\beta_p=p^{2k-1}$.
Explicitly one has
\begin{align}
\label{w3.4}
\alpha_p,\beta_p=\frac{\lambda_p\pm \sqrt{\lambda_p^2-4p^{2k-1}}}{2}.
\end{align}
Now assume that $a(tp^{2m})\geq 0$ for almost all $m$. Then by Landau's theorem 
the Dirichlet series on the left hand side of \eqref{w3.3} 
either converges everywhere or 
has a singularity at the real point of its abscissa of convergence. 
The first alternative clearly is impossible, since the right-hand side of 
\eqref{w3.3} has a pole for $p^s=\alpha_p$ or $p^s=\beta_p$.

Thus the second alternative must hold, and in particular 
$\alpha_p$ or $\beta_p$ must be real. By Deligne's theorem, we have
\[
\lambda_p^2\leq 4p^{2k-1},
\]
hence in combination with \eqref{w3.4} we find that
\[
\lambda_p=\pm 2p^{k-1/2}.
\]
In particular we conclude that $\sqrt p$ is contained in the number field 
$K_f$ generated by the Hecke eigenvalues of $f$.

Since for different primes $p_1,\dots ,p_r$ the degree of the field 
extension 
\[
\Q(\sqrt{p_1},\dots ,\sqrt{p_r})/\Q
\]
is $2^r$, 
we deduce our assertion.
\end{proof}

Throughout the rest of this section, let $f=\sum_{n=1}^\infty a(n)
q^n\in S_{k+1/2}(N,\chi)$ be an arbitrary  
non-zero cup form with Fourier
coefficients $a(n)\in \R$.
For the proof of Theorem \ref{thm1} we need the following three
propositions.

\begin{proposition}
\label{prop0} There exist infinitely many $n\in \N$ such that $a(n)$
is negative.
\end{proposition}

\begin{proof}
  This result is proved in \cite{KKP} in much greater generality. For
  the convenience of the reader we sketch the argument in the present
  special case.

Assume that there exist only finitely many $n\in \N$ such that
$a(n)<0$. Then the Hecke $L$-function
of $f$ (which is entire, cf.~\cite{Sh}) converges for all $s\in \C$
by Landau's theorem. Consequently,
\[
a(n) \ll_C n^C
\]
for all $C\in \R$.

This implies that the Rankin $L$-function of $f$ also converges for
all $s\in \C$. Arguing as at the end of the proof of Theorem
\ref{thm2} we find that $f$ vanishes identically, contradicting the
assumption $f\neq 0$.
%
\end{proof}

\begin{proposition}
\label{prop1} Let $p$ be a prime not dividing $N$, and let $\eps\in
\{\pm 1\}$. Assume that $a(n)\geq 0$ for all positive integers $n$
with $(\frac{n}{p})=\eps$. Then $f$ is an eigenform of the Hecke
operator $T(p^2)$ with eigenvalue $-\eps \chi^*(p)(p^k+p^{k-1})$.
\end{proposition}

\begin{proof}
We consider the  cusp form
\[
\tilde f
:= \sum_{\substack{n\geq 1\\ \leg{n}{p} =\eps}} a(n) q^n.
\]
According to \cite{Br}, Section 2 (iii) and (iv), it belongs to the
space $S_{k+1/2}(N p^2,\chi)$. By assumption, $\tilde f$ has
non-negative real coefficients. Using Proposition \ref{prop0} we
find that $\tilde f=0$. Consequently, $a(n)=0$ for all positive
integers $n$ with $(\frac{n}{p})=\eps$. Now the assertion follows
from \cite{Br}, Lemma 1.
\end{proof}

\begin{proposition}
\label{prop2} Let $p$ be a prime not dividing $N$, and let $\eps\in
\{\pm 1\}$. Assume that $f$ is contained in $S_{k+1/2}^*(N,\chi)$.
Then there exist positive integers $n,n'$ such that
\begin{align*}
\leg{n}{p}&=\eps\quad\text{and $a(n)<0$,}\\
\leg{n'}{p}&=\eps\quad\text{and $a(n')>0$.}
\end{align*}
\end{proposition}

\begin{proof}
Suppose that $a(n)\geq 0$ for all $n\in \N$ with
$(\frac{n}{p})=\eps$.
%
Then, by Proposition \ref{prop1}, $f$ is an eigenform of $T(p^2)$
with eigenvalue $\lambda_p=-\eps \chi^*(p)(p^k+p^{k-1})$. Using the
Shimura lift, we see that $\lambda_p$ is also an eigenvalue of the
integral weight Hecke operator $T(p)$ on $S_{2k}(N/2,\chi^2)$. But
it is easy to see that any eigenvalue $\lambda$ of this Hecke
operator satisfies the bound $|\lambda|<p^k+p^{k-1}$, see
e.g.~\cite{Ko} (or use the stronger Deligne bound). We obtain a
contradiction.

Finally, replacing $f$ by $-f$ we deduce the existence of $n'$ with
the claimed properties.
\end{proof}

\begin{proof}[Proof of Theorem \ref{thm1}.]
Suppose that there exist finitely many square-free
$t_1,\dots,t_h\in\N$ such that $a(tn^2)\leq 0$ for all square-free
integers $t$ different from $t_\nu$, $\nu=1,\dots,h$, and all
$n\in\N$. Choose a prime $p$ coprime to $N$ such that
\[
\leg{t_\nu}{p}= 1,\quad \text{for all $\nu=1,\dots,h$.}
\]
Then $a(n)\leq 0 $ for all $n\in \N$ with $(\frac{n}{p})=-1$. But
this contradicts Proposition \ref{prop2}. Hence there exist
infinitely many square-free $t\in \N$ for which there is an $n_t\in
\N$ such that $a(t n_t^2)>0$.

Finally, replacing $f$ by $-f$ we deduce the existence of infinitely
many square-free $t\in \N$ for which there is an $n_t\in \N$ such
that $a(t n_t^2)<0$.
\end{proof}

\begin{proof}[Proof of Theorem \ref{cor1}.]
The assertion follows combining Theorem \ref{thm1} and
Proposition~\ref{prop2}.
\end{proof}

\section{Examples}

\medskip

Let $f\in S_{k+1/2}^*(N,\chi_0)$ be a cusp form with trivial
character $\chi_0$, square-free level, and real coefficients $a(n)$.
We suppose that $f$ is contained in the plus space, that is,
$a(n)=0$ when $(-1)^k n\equiv 2,3\pmod{4}$, see \cite{KZ},
\cite{Ko2}. For a positive number $X$, we define the quantity
\[
R^{+}_{tot}(f,X)=\frac{\#\{n\leq X;\quad a(n) >0\}}{\#\{n\leq
  X;\quad a(n) \neq 0\}}.
\]
Moreover, in view of Waldspurger's theorem, it is natural to
consider the coefficients $a(d)$ especially for {\em fundamental}
discriminants $d$.  Therefore we put
\[
R^{+}_{fund}(f,X)=\frac{\#\{d\leq X;\quad \text{$d$ fundamental
    discriminant and $a(d) >0$}\}}{\#\{d\leq X;\quad \text{$d$
    fundamental discriminant and $a(d) \neq 0$}\}}.
\]
The numerical experiments below suggest that
\[
\lim_{X\to
  \infty}R^{+}_{tot}(f,X)=1/2,\qquad \lim_{X\to
  \infty}R^{+}_{fund}(f,X)=1/2.
\]

In our first example, we consider the Delta function
$\Delta(z)=q\prod_{n\geq 1} (1-q^n)^{24}$. A cusp form of weight
$13/2$ of level $4$ in the plus space corresponding to $\Delta$
under the Shimura lift is given by
\[
\delta(z)=\frac{1}{8\pi i} \big(2 E_4(4z)\theta'(z)
-E'_4(4z)\theta(z)\big)\in S_{13/2}^+(4,\chi_0),
\]
see \cite{KZ}. Here $E_4(z)=1+240\sum_{n\geq 1} \sigma_3(n)q^n$ is
the classical Eisenstein series of weight $4$ and
$\theta(z)=\sum_{n\in \Z} q^{n^2}$. The Fourier expansion of
$\delta$ starts as follows:
\[
\delta(z)= q -56 q^4+ 120 q^5 - 240 q^8 +9 q^9 + 1440 q ^{12} -1320
q^{13}-704 q^{16}-240 q^{17}+ \dots .
\]
Computational data for $\delta$ is listed in Table \ref{table1}.

\begin{table}[h]
\caption{\label{table1} The proportion of positive coefficients of
$\delta$}
\begin{tabular}{|r||c|c|c|c|c|c| }
\hline \rule[-3mm]{0mm}{8mm}
$X$ & 10 & $10^2$ &   $10^3$ & $10^4$ & $10^5$ &  $10^6$ \\
\hline \rule[-3mm]{0mm}{8mm} $R^{+}_{tot}(\delta,X)$ & $0.600$ &
$0.520$ & $0.518$ & $0.504600$ & $0.499600$ &
0.499822   \\
\hline \rule[-3mm]{0mm}{8mm} $R^{+}_{fund}(\delta,X)$ & $0.667$ &
$0.548$ & $0.515$ & $0.501643$ & $0.500016$
 & 0.499836   \\
\hline
\end{tabular}
\end{table}

In our second example, we consider the cusp form $G=\eta(z)^2\eta(11
z)^2$ of weight $2$ and level $11$ corresponding to the elliptic
curve $X_0(11)$. Here $\eta=q^{1/24}\prod_{n\geq 1} (1-q^n)$ is the
Dedekind eta function.  A cusp form of weight $3/2$ and level $44$
in the  plus space corresponding to $G$ under the Shimura lift is
\[
g(z) = \big(\theta(11z)\eta(2z)\eta(22z)\big)| U_4 \in
S^+_{3/2}(44,\chi_0),
\]
see \cite{Du} \S2.  Here  $U_4$ denotes the usual Hecke operator of
index $4$.  The Fourier expansion of $g$ starts as follows:
\[
g(z)=
q^3-q^4-q^{11}-q^{12}+q^{15}+2q^{16}+q^{20}-q^{23}-q^{27}-q^{31}+q^{44}+q^{55}
+\dots.
\]
Computational data for $g$ is listed in Table \ref{table2}.

\begin{table}[h]
\caption{\label{table2} The proportion of positive coefficients of
$g$}
\begin{tabular}{|r||c|c|c|c|c|c|c| }
\hline \rule[-3mm]{0mm}{8mm}
$X$ & 10 & $10^2$ &   $10^3$ & $10^4$ & $10^5$ &  $10^6$ \\
\hline \rule[-3mm]{0mm}{8mm}
$R^{+}_{tot}(g,X)$  & 0.500 & 0.500 & 0.500 & $0.496042$ & $0.501022$ & $0.499544$    \\
\hline \rule[-3mm]{0mm}{8mm}
$R^{+}_{fund}(g,X)$ & 1.000 & 0.500 & 0.503 & $0.491968$ & $0.500861$ & $0.499589$    \\
\hline
\end{tabular}
\end{table}

\end{document}